\documentclass[12pt]{amsart}
\usepackage{amsthm,amsfonts,amssymb,amscd,mathrsfs}
\usepackage{bezier,longtable,amstext}
\usepackage{amssymb}

\usepackage[all]{xy}

\newcommand\CC{{\mathbb{C}}}
\newcommand\PP{{\mathbb{P}}}

\newcommand\coker{\operatorname{coker}\nolimits}

\newcommand\Hom{\operatorname{Hom}\nolimits}

\newcommand\im{\operatorname{Im}\nolimits}

\renewcommand\emptyset{\varnothing}

\newcommand\lra{\longrightarrow}
\newcommand\lrdash{\:
\xymatrix@1{\ar@{-->}[r]&}\:
}
\newcommand{\lrdashar}[1]{\:
\xymatrix@1{\ar@{-->}[r]^{#1}&}\:
}
\newcommand{\sendsto}[1]{\:
\xymatrix@1{\ar@{|->}[r]^{#1}&}\:
}
\newcommand{\intoo}[1]{\:
\xymatrix@1{\ar@{^(->}[r]^{#1}&}\:
}

\newtheorem{theorem}{Theorem}[section]
\newtheorem{proposition}[theorem]{Proposition}

\newtheorem{lemma}[theorem]{Lemma}

\newtheorem{sub}[subsection]{}

\theoremstyle{definition}
\newtheorem{definition}[theorem]{Definition}

\newtheorem{proposition-definition}[theorem]{Proposition-Definition}
\newtheorem{remark}[theorem]{Remark}

\nonstopmode

\begin{document}

\vspace{1cm}

\title[Symplectic instantons on $\mathbb{P}^3$]{Moduli of symplectic instanton vector bundles\\[5pt] of higher rank on
projective space $\mathbb{P}^3$}

\author{U. Bruzzo}
\address{\scriptsize U. Bruzzo:
Scuola Internazionale Superiore di Studi Avanzati,
Via Bonomea 265, 34136 Trieste, Italia  and  Istituto Nazionale di Fisica Nucleare, sezione di Trieste}
\email{bruzzo@sissa.it}

\author{D. Markushevich}
\address{\scriptsize D. Markushevich:
Math\'ematiques - b\^{a}t.M2, Universit\'e Lille 1,
F-59655 Villeneuve d'Ascq Cedex, France}
\email{markushe@math.univ-lille1.fr}

\author{A.~S. Tikhomirov}
\address{\scriptsize   A.S.~Tikhomirov:
Department of Mathematics\\
Yaroslavl State Pedagogical University\\
Respublikanskaya Str. 108\\
150 000 Yaroslavl, Russia}
\email{astikhomirov@mail.ru}

\thispagestyle{empty}

\begin{abstract}

Symplectic instanton vector bundles on the projective space $\mathbb{P}^3$ constitute a natural generalization
of mathematical instantons of rank 2. We study the moduli space
$I_{n,r}$ of rank-$2r$ symplectic instanton vector bundles on $\mathbb{P}^3$ with $r\ge2$ and second Chern class
$n\ge r,\ n\equiv r({\rm mod}2)$.
We give an explicit construction of an irreducible component $I^*_{n,r}$ of this space for each such value of $n$
and show that $I^*_{n,r}$ has the expected dimension $4n(r+1)-r(2r+1)$.

\end{abstract}

\maketitle

\section{Introduction}\label{sec0}

By a {\it symplectic instanton vector bundle} of rank $2r$ and charge $n$ (shortly, a {\it symplectic
$(n,r)$-instanton}) on the 3-dimensional projective space $\mathbb{P}^3$ we understand
an algebraic vector bundle $E=E_{2r}$  of rank $2r$ on $\mathbb{P}^3$ with Chern classes
\begin{equation}\label{1st Chern class}
c_1(E)=0,
\end{equation}
\begin{equation}\label{2nd Chern class}
c_2(E)=n,\ \ \ n\ge1,
\end{equation}
supplied with a symplectic structure and satisfying the vanishing conditions
\begin{equation}\label{vanishing condns}
h^0(E)=h^1(E\otimes\mathcal{O}_{\mathbb{P}^3}(-2))=0.
\end{equation}
By a symplectic structure we mean an anti-self-dual
isomorphism
\begin{equation}\label{sympl str}
\phi:\ E\overset{\simeq}\to E^\vee,\ \ \ \phi^\vee=-\phi,
\end{equation}
considered modulo proportionality. The vanishing of the first Chern class \eqref{1st Chern class} follows from the
existence of a symplectic structure (\ref{sympl str}), and if $r=1$, then the two conditions are equivalent. We will
denote the moduli space of symplectic $(n,r)$-instantons  by $I_{n,r}$.

For $r=1$ these bundles  relate, via the so-called Atiyah-Ward correspondence, to rank-2 ``physical"  instantons over
the 4-sphere $S^4$, these   being  anti-self-dual connections with structure group $SU(2)=\mathbf{Sp}(1)$ \cite{AW}.
Important results on the moduli spaces $I_{n}=I_{n,1}$ of rank-2 instantons have been obtained recently: smoothness
\cite{JV} for all $n$, irreducibility \cite{T} for odd $n$.

Much less is known about the moduli spaces $I_{n,r}$ for $r>1$. In fact the symplectic instantons with $r>1$ are
as natural as those with $r=1$, for they are related, via the same Atiyah-Ward correspondence, to the anti-self-dual
connections over $S^4$ with structure group $\mathbf{Sp}(r)$, see \cite{A}. As far as we know, the present paper is
the first one addressing the properties of the corresponding spaces $I_{n,r}$. The tool we use to construct $I_{n,r}$
is the monad method; it originates in the work of Horrocks \cite{H} and is known as the ADHM construction of instantons
since \cite{ADHM}. It was further sharpened in the work of Barth \cite{B},
Barth and Hulek \cite{BH} and Tyurin \cite{Tju1}, \cite{Tju2}. This method permits to encode the instantons, usual ones
or symplectic of higher rank, by hyperwebs of quadrics.

For a sample of the physical literature about symplectic instantons, see e.g.~\cite{Mc}.

We fix basic terminology and notation in Section 2 and introduce the hyperwebs of quadrics in Section 3.
We prove that, for any $r\ge2$ and for any $n\ge r$ such that $n\equiv r({\rm mod\:}2),$ the moduli space
$I_{n,r}$ is nonempty and is realized as a free quotient $MI_{n,r}/(GL(n)/\pm{\rm id})$, where $MI_{n,r}$ is a
Zariski locally closed subset of an affine space (see Theorem \ref{princbundle}). Thus $MI_{n,r}$ carries a natural
structure of a reduced scheme, and $I_{n,r}$ is an algebraic space.
In Section 4 we give an explicit construction of vector bundles from $I_{n,r}$ for the above values of $n$ and $r$ and
introduce a component $I^*_{n,r}$ of $I_{n,r}$ characterized by a certain open condition
(*), see Definition \ref{property (*)}.
In Section \ref{I^*nr} we prove Theorem \ref{Irred of I*} on the irreducibility of
$I^*_{n,r}$,  the main result of this paper.

\medskip

\textbf{Acknowledgements.}   D.~M.~was partially supported by the grant VHSMOD-2009 No. ANR-09-BLAN-0104-01,
and U.~B.~by PRIN ``Geometria delle variet\`a algebriche e dei loro spazi di moduli". U.~B.~and A.~S.~T.~acknowledge
support and hospitality of the Max Planck Institute for Mathematics in Bonn, where they started the work on this paper
during their stay in winter 2011. U.~B.~is a member of the VBAC group.

\bigskip
\section{Notation and conventions}\label{notation}

In many respects, we follow the exposition of \cite{T}, and we stick to the notation introduced therein.
The base field $\mathbf{k}$ is assumed to be
algebraically closed of characteristic 0. We identify vector bundles with locally free
sheaves. If $\mathcal{F}$ is a sheaf of $\mathcal{O}_X$-modules on an algebraic
variety or a scheme $X$, then $n\mathcal{F}$ denotes a direct sum of $n$ copies of $\mathcal{F}$, $H^i(\mathcal{F})$
denotes the $i^{th}$ cohomology group of
$\mathcal{F}$, $h^i(\mathcal{F}):=\dim H^i(\mathcal{F})$, and $\mathcal{F}^\vee$
denotes the dual of $\mathcal{F}$, that is, $\mathcal{F}^\vee:=
\mathcal{H}om_{\mathcal{O}_X} (\mathcal{F},\mathcal{O}_X)$. If
$X=\mathbb{P}^r$ and $t$ is an integer, then by $\mathcal{F}(t)$ we denote the sheaf
$\mathcal{F}\otimes\mathcal{O}_{\mathbb{P}^r}(t)$. $[\mathcal{F}]$ will denote the
isomorphism class of a sheaf $\mathcal{F}$. For any morphism of
$\mathcal{O}_X$-sheaves $f:\mathcal{F}\to\mathcal{F}'$ and any $\mathbf{k}$-vector space $U$
(respectively, for any homomorphism $f:U\to U'$ of $\mathbf{k}$-vector spaces) we will denote,
for short, by the same letter $f$ the induced morphism of sheaves
$id\otimes f:U\otimes\mathcal{F}\to U\otimes\mathcal{F}'$
(respectively, the induced morphism
$f\otimes id:U\otimes\mathcal{F}\to U'\otimes\mathcal{F}$).

We fix an integer $n\ge1$ and denote by $H_n$ a fixed $n$-dimensional vector space over $\mathbf{k}$.
Throughout the paper, $V$ will be a fixed vector space of dimension 4 over $\mathbf{k}$, and
we set $\mathbb{P}^3:=P(V)$. We reserve the letters $u$ and $v$ for
denoting the two morphisms in the Euler exact sequence
$0\to\mathcal{O}_{\mathbb{P}^3}(-1)\overset{u}\to V^\vee\otimes\mathcal{O}_{\mathbb{P}^3}
\overset{v}\to T_{\mathbb{P}^3}(-1)\to0$.
For any $\mathbf{k}$-vector spaces $U$ and $W$ and any vector
$\phi\in{\rm Hom}(U,W\otimes\wedge^2V^\vee)\subset{\rm Hom}(U\otimes V,W\otimes V^\vee)$
understood as a linear map $\phi:U\otimes V\to W\otimes V^\vee$ or, equivalently,
as a map
${}^\sharp\phi:U\to W\otimes\wedge^2V^\vee$, we will denote by $\widetilde{\phi}$ the composition
$U\otimes\mathcal{O}_{\mathbb{P}^3}\overset{{}^\sharp\phi}\to
W\otimes\wedge^2V^\vee\otimes\mathcal{O}_{\mathbb{P}^3}\overset{\epsilon}\to
W\otimes\Omega_{\mathbb{P}^3}(2)$, where $\epsilon$ is the induced morphism in the exact triple
$0\to\wedge^2\Omega_{\mathbb{P}^3}(2)\overset{\wedge^2v^\vee}\to\wedge^2V^\vee\otimes
\mathcal{O}_{\mathbb{P}^3}\overset{\epsilon}\to\Omega_{\mathbb{P}^3}(2)\to0$ obtained by taking
the second wedge power of the dual Euler exact sequence.

Given an integer $m\ge1$, we denote by $\mathbf{S}_m$ (resp. $\mathbf{\Sigma}_{m+1}$)  the vector space
$S^2H_m^\vee\otimes \wedge^2V^\vee$ (resp.
$\Hom(H_m,H_{m+1}^\vee\otimes\wedge^2V^\vee)$). Abusing notation, we will denote by the same symbol a
$\mathbf{k}$-vector space, say $U$, and the associated affine space $\mathbf{V}(U^\vee)={\rm Spec}(Sym^*U^\vee)$.

All the schemes considered in the paper are Noetherian. By a general point of an irreducible
(but not necessarily reduced) scheme $\mathcal{X}$ we mean any closed point of some dense open subset of
$\mathcal{X}$. An irreducible scheme is called generically reduced if it is reduced at a general point.

\bigskip

\section{Generalities on symplectic instantons and definition of $MI_{n,r}$}\label{general}

In this section we enumerate some facts about symplectic instantons which are completely parallel to those for
rank-2 usual instantons, see \cite[Section 3]{T}.

For a given symplectic $(n,r)$-instanton $E$, the first condition (\ref{vanishing condns}) yields
$h^0(E(-i))=0, i\ge0,$ which, together with the exact triple
$0\to E(-j-1)\to E(-j)\to E(-j)|_{\mathbb{P}^2}\to0$
for $j=0$ and (\ref{vanishing condns}), implies that $h^0(E(-1)|_{\mathbb{P}^2})=0$, hence also
$h^0(E(-i)|_{\mathbb{P}^2})=0,\ i\ge1$. The last equality for
$i=2$, together with (\ref{vanishing condns}) and the above triple for $j=2$, gives $h^1(E(-3))=0$,
hence also $h^1(E(-4))=0$. Then, from Serre duality and (\ref{sympl str}), we deduce:
\begin{equation}\label{dimension2}
h^i(E)=h^i(E(-1))=h^{3-i}(E(-3))=h^{3-i}(E(-4))=0,\ \ i\ne1,\ \ \ h^i(E(-2))=0,\ \ i\ge0.
\end{equation}
By Riemann-Roch and (\ref{vanishing condns}), (\ref{dimension2}), we have
\begin{equation}\label{dimension1,5}
h^1(E(-1))=h^2(E(-3))=n,\ h^1(E)=h^2(E(-4))=2n-2r. \end{equation}
By tensoring the dual Euler sequence by $E$ we also obtain
\begin{equation}\label{dimension1}
h^1(E\otimes\Omega^1_{\mathbb{P}^3})=h^2(E\otimes\Omega^2_{\mathbb{P}^3})=2n+2r,\ \ \
\end{equation}

Consider a triple $(E,f,\phi)$ where $E$ is a $(n,r)$-instanton, $f:
H_n\overset{\simeq}\to H^2(E(-3))$ an isomorphism
and $\phi:E\overset{\simeq}\to E^\vee$ a symplectic structure on $E$.
Two triples
$(E,f,\phi)$ and $(E'f',\phi')$ are called equivalent if there is an isomorphism
$g:E\overset{\simeq}\to E'$ such that $g_*\circ f=\lambda f'$ with $\lambda\in\{1,-1\}$ and
$\phi=g^\vee\circ \phi'\circ g$,
where
$g_*:H^2(E(-3))\overset{\simeq}\to H^2(E'(-3))$
is the induced isomorphism. We denote by $[E,f,\phi]$ the equivalence class of a triple $(E,f,\phi)$.
It follows from this definition that the set $F_{[E]}$ of all equivalence classes
$[E,f,\phi]$ with given $[E]$ is a homogeneous space of the group $GL(H_n)/\{\pm{\rm id}\}$.

Each class $[E,f,\phi]$ defines a point
\begin{equation}\label{hypernet}
A=A([E,f,\phi])\in S^2H_n^\vee\otimes \wedge^2V^\vee
\end{equation}
in the following way. Consider the exact sequences
\begin{equation}\label{Koszul triples}
0\to\Omega^1_{\mathbb{P}^3}\overset{i_1}\to V^\vee\otimes\mathcal{O}_{\mathbb{P}^3}(-1)\to
\mathcal{O}_{\mathbb{P}^3}\to0,
\end{equation}
$$
0\to\Omega^2_{\mathbb{P}^3}\to\wedge^2V^\vee\otimes\mathcal{O}_{\mathbb{P}^3}(-2)\to
\Omega^1_{\mathbb{P}^3}\to0,
$$
$$
0\to\wedge^4V^\vee\otimes\mathcal{O}_{\mathbb{P}^3}(-4)\to
\wedge^3V^\vee\otimes\mathcal{O}_{\mathbb{P}^3}(-3)\overset{i_2}\to\Omega^2_{\mathbb{P}^3}\to0,
$$
induced by the Koszul complex of
$V^\vee\otimes\mathcal{O}_{\mathbb{P}^3}(-1)\overset{ev}
\twoheadrightarrow\mathcal{O}_{\mathbb{P}^3}$.
Twisting these sequences by $E$ and taking into account (\ref{vanishing condns}), \eqref{dimension2}-\eqref{dimension1}, we obtain the vanishing
\begin{equation}\label{dimension3}
h^0(E\otimes\Omega_{\mathbb{P}^3})=h^3(E\otimes\Omega^2_{\mathbb{P}^3})=
h^2(E\otimes\Omega_{\mathbb{P}^3})=0
\end{equation}
and the diagram with exact rows
\begin{equation}\label{A'}
\xymatrix{0\ar[r] &
H^2(E(-4))\otimes\wedge^4V^\vee\ar[r] &
H^2(E(-3))\otimes\wedge^3V^\vee\ar[r]^{\ \ \ \ i_2}\ar[d]^{A'} &
H^2(E\otimes\Omega^2_{\mathbb{P}^3})\ar[r]& 0 \\
0 & H^1(E))\ar[l] & H^1(E(-1))\otimes V^\vee\ar[l] &
H^1(E\otimes\Omega_{\mathbb{P}^3})\ar[l]_{\ \ \ \ i_1}\ar[u]^{\cong}_{\partial}& 0,\ar[l]}
\end{equation}
where $A':=i_1\circ\partial^{-1}\circ i_2$. The Euler exact sequence (\ref{Koszul triples})
yields the canonical isomorphism
$\omega_{\mathbb{P}^3}\overset{\simeq}\to\wedge^4V^\vee\otimes\mathcal{O}_{\mathbb{P}^3}(-4)$,
and fixing an isomorphism
$\tau:\mathbf{k}\overset{\simeq}\to\wedge^4V^\vee$
we have the isomorphisms
$\tilde{\tau}:V\overset{\simeq}\to\wedge^3V^\vee$
and
$\hat{\tau}:\omega_{\mathbb{P}^3}\overset{\simeq}\to\mathcal{O}_{\mathbb{P}^3}(-4)$.
We define $A$ in (\ref{hypernet}) as the composition
\begin{equation}\label{An}
A:H_n\otimes V\overset{\tilde{\tau}}{\overset{\simeq}\to}
H_n\otimes\wedge^3V^\vee\overset{f}{\overset{\simeq}\to}
H^2(E(-3))\otimes\wedge^3V^\vee\overset{A'}\to H^1(E(-1))\otimes V^\vee
\overset{\phi}{\overset{\simeq}\to}
\end{equation}
$$
\overset{\phi}{\overset{\simeq}\to} H^1(E^\vee(-1))\otimes V^\vee
\overset{SD}{\overset{\simeq}\to} H^2(E(1)\otimes\omega_{\mathbb{P}^3})^\vee\otimes V^\vee
\overset{\hat{\tau}}{\overset{\simeq}\to} H^2(E(-3))^\vee\otimes V^\vee
\overset{f^\vee}{\overset{\simeq}\to}H_n^\vee\otimes V^\vee,
$$
where $SD$ is the Serre duality isomorphism. One can verify that $A$ is a skew symmetric map
which depends only on the class $[E,f,\phi]$, but does not depend on the choice of $\tau$, and that $A\in\wedge^2(H_n^\vee\otimes V^\vee)$ lies in the direct summand
$\mathbf{S}_n=S^2H_n^\vee\otimes \wedge^2V^\vee$
of the canonical decomposition
\begin{equation}\label{can decomp}
\wedge^2(H_n^\vee\otimes V^\vee)=S^2H_n^\vee\otimes \wedge^2V^\vee\oplus
\wedge^2H_n^\vee\otimes S^2V^\vee.
\end{equation}
Here $\mathbf{S}_n$ is the space of hyperwebs of quadrics in $H_n$. For this reason we call $A$ the
$(n,r)$-{\it instanton hyperweb of quadrics} corresponding to the data $[E,f,\phi]$.

Denote $W_A:=H_n\otimes V/\ker A$. Using the above chain of isomorphisms we
can rewrite the diagram (\ref{A'}) as
\begin{equation}\label{qA}
\xymatrix{0\ar[r] & \ker A\ar[r] &
H_n\otimes V\ar[r]^{\ \ c_A}\ar[d]^{A} & W_{A}\ar[r]\ar[d]_{\cong}^{q_{A}} & 0 \\
0 & \ker A^\vee\ar[l] & H_n^\vee\otimes V^\vee\ar[l] &
W_{A}^\vee\ar[l]_{\ \ \ \ \ c_A^\vee}& 0.\ar[l]}
\end{equation}
In view of (\ref{dimension1}), $\dim W_{A}=2n+2r$ and
$q_{A}:W_{A}\overset{\simeq}\to W_{A}^\vee$
is a skew-symmetric isomorphism. An important property of $A=A([E,f,\phi])$ is that
the induced morphism of sheaves
\begin{equation}\label{an^vee}
a_A^\vee:{W}^\vee_{A}\otimes\mathcal{O}_{\mathbb{P}^3}
\overset{c_A^\vee}\to
H_n^\vee\otimes V^\vee\otimes\mathcal{O}_{\mathbb{P}^3}
\overset{ev}\to H_n^\vee\otimes\mathcal{O}_{\mathbb{P}^3}(1)
\end{equation}
is surjective and the composition
$H_n\otimes\mathcal{O}_{\mathbb{P}^3}(-1)\overset{a_A}\to
W_{A}\otimes\mathcal{O}_{\mathbb{P}^3}\overset{q_A}\to
W^\vee_{A}\otimes\mathcal{O}_{\mathbb{P}^3}
\overset{a^\vee_{A}}\to H_n^\vee\otimes\mathcal{O}_{\mathbb{P}^3}(1)$
is zero. Applying Beilinson spectral sequence \cite{Bei} to $E(-1)$, we see that $E\simeq\ker(a_A^\vee\circ q_{A})/\im a_A$. Thus $A$ defines a monad
\begin{equation}\label{Monad A}
\mathcal{M}_{A}:\ \ 0\to H_n\otimes\mathcal{O}_{\mathbb{P}^3}(-1)\overset{a_A}\to
W_{A}\otimes\mathcal{O}_{\mathbb{P}^3}\overset{a_A^\vee\circ q_{A}}
\to H_n^\vee\otimes\mathcal{O}_{\mathbb{P}^3}(1)\to0\ ,
\end{equation}
whose cohomology sheaf
\begin{equation}\label{coho sheaf}
E_{2r}(A):=\ker(a_A^\vee\circ q_{A})/\im a_A.
\end{equation}
is isomorphic to $E$.
Twisting $\mathcal{M}_{A}$ by $\mathcal{O}_{\mathbb{P}^3}(-3)$ and using (\ref{coho sheaf}), we obtain the isomorphism
$f:H_n\overset{\simeq}\to H^2(E(-3))$.
Furthermore, the fact that $q_{A}$ is symplectic implies
that there is a canonical isomorphism of $\mathcal{M}_{A}$ with its dual which induces the
symplectic isomorphism $\phi:E\overset{\simeq}\to E^\vee$. Thus, the data $[E,f,\phi]$ are recovered
from $A$. This leads to the following description of the moduli space $I_{n,r}$. Consider
the {\it set of $(n,r)$-instanton hyperwebs of quadrics}
\renewcommand\theenumi{\roman{enumi}}
\begin{equation}\label{space of nets}
MI_{n,r}:=\left\{A\in \mathbf{S}_n\ \left|\
\begin{minipage}{26em}
\begin{enumerate}
\item $rk(A:H_n\otimes V\to H_n^\vee\otimes V^\vee)=2n+2r$,
\item the morphism $a_A^\vee:W^\vee_{A}\otimes\mathcal{O}_{\mathbb{P}^3}\to
H_n^\vee\otimes\mathcal{O}_{\mathbb{P}^3}(1)$ defined by $A$ in (\ref{an^vee}) is surjective,
\item $h^0(E_{2r}(A))=0$, where $E_{2r}(A)=
\ker(a^\vee_A\circ q_{A})/\im a_A$ and $q_{A}:W_{A}\overset{\simeq}\to W_{A}^\vee$ is a symplectic
isomorphism associated to $A$ by (\ref{qA}).
\end{enumerate}
\end{minipage}\right.
\right\}
\end{equation}
It is a locally closed subscheme
of the affine space $\mathbf{S}_n$.

\begin{theorem}\label{princbundle}
The natural morphism
\begin{equation}\label{pi nr}
\pi_{n,r}: MI_{n,r}\to I_{n,r},\ A\mapsto[E_{2r}(A)],
\end{equation}
is a principal
$GL(H_n)/\{\pm{\rm id}\}$-bundle in the \'etale topology.
Hence $I_{n,r}$ is a  quotient stack $MI_{n,r}/(GL(H_n)/\{\pm{\rm id}\})$, making it an algebraic space.\end{theorem}
\begin{proof} See \cite[Section 3]{T}.\end{proof}

Each fibre $F_{[E]}=\pi_n^{-1}([E])$ over an arbitrary point $[E]\in I_{n,r}$ is a principal homogeneous
space of the group $GL(H_n)/\{\pm{\rm id}\}$. Hence the
irreducibility of $(I_{n,r})_{red}$ is equivalent to the irreducibility of the scheme $(MI_{n,r})_{red}$.

We can also state:
\begin{theorem}\label{rank condns}
For each $n\ge1$, the space $MI_{n,r}$ of $(n,r)$-instanton hyperwebs of quadrics is a locally closed subscheme
of the vector space $\mathbf{S}_n$ given locally at any point
$A\in MI_{n,r}$ by
\begin{equation}\label{eqns(i)}
\binom{2n-2r}{2}=2n^2-n(4r+1)+r(2r+1)
\end{equation}
equations obtained as the rank condition (i) in (\ref{space of nets}).
\end{theorem}

Note that from (\ref{eqns(i)}) it follows that
\begin{equation}\label{dim MInr ge...}
\dim_{[A]}MI_{n,r}\ge\dim \mathbf{S}_n-(2n^2-n(4r+1)+r(2r+1))=n^2+4n(r+1)-r(2r+1)
\end{equation}
at any point $A\in MI_{n,r}$. Hence,
\begin{equation}\label{dim Inr ge...}
\dim_{[E]}I_{n,r}\ge4n(r+1)-r(2r+1)
\end{equation}
at any point $[E]\in I_{n,r}$, since $MI_{n,r}\to I_{n,r}$ is a principal
$GL(H_n)/\{\pm{\rm id}\}$-bundle in the \'etale topology.

\bigskip

\section{Explicit construction of symplectic instantons}\label{explicit}

\begin{sub}{\bf Example: symplectic $(n,n)$-instantons.}
\label{nn-instantons}
\rm
In this subsection we recall some known facts about symplectic $(n,n)$-instantons and their relation to usual rank-2
instantons, see \cite[Sections 5-6]{T}.
We first show that each invertible hyperweb of quadrics $A\in\mathbf{S}_n$ naturally leads to a construction
of a symplectic $(n,n)$-instanton $E_{2n}(A)$ on $\mathbb{P}^3$.
Given an integer $n\ge1$, set
\begin{equation}\label{S0n}
\mathbf{S}^0_n:=\{A\in\mathbf{S}_n\ |\
A:H_n\otimes V\to H_n^\vee\otimes V^\vee
\ {\rm is\ an\ invertible\ map}\}.
\end{equation}
Then $\mathbf{S}^0_n$ is a dense open subset of $\mathbf{S}_n$,
and it is easy to see that for any $A\in\mathbf{S}^0_n$ the following conditions are satisfied.

(1) The morphism
$\widetilde{A}:H_n\otimes\mathcal{O}_{\mathbb{P}^3}(-1)\to
H_n^\vee\otimes\Omega_{\mathbb{P}^3}(1)$
induced by
$A$
is a subbundle embedding, and
\begin{equation}\label{E2n}
E_{2n}(A):=\coker(\widetilde{A})
\end{equation}
is a symplectic $(n,n)$-instanton, that is,
\begin{equation}\label{E2n in Inn}
[E_{2n}(A)]\in I_{n,n}.
\end{equation}

(2) For all $i\geq 0$,
\begin{equation}\label{vanish cohom of E2n}
h^i(E_{2n}(A))=h^i(E_{2n}(A)(-2))=0.
\end{equation}

This follows from the diagram
\begin{equation}\label{diag E2n}
\xymatrix{
& & 0\ar[d] & 0\ar[d] & & \\
& 0\ar[r] & H_n\otimes\mathcal{O}_{\mathbb{P}^3}(-1)
\ar[r]^{\widetilde{A}}\ar[d]^u
& H_n^\vee\otimes\Omega_{\mathbb{P}^3}(1)\ar[r]^{\ \ \ e}\ar[d]^{v^\vee} &
E_{2n}(A)\ar[r] &  0\\
& &H_n\otimes V\otimes\mathcal{O}_{\mathbb{P}^3}
\ar[r]^{A\ }_-{\simeq}\ar[d]^v
& H_n^\vee\otimes V^\vee\otimes\mathcal{O}_{\mathbb{P}^3}\ar[d]^{u^\vee} &  & \\
&0\to E_{2n}(A)^\vee\ar[r] &
H_n\otimes T_{\mathbb{P}^3}(-1)\ar[r]^-{\widetilde{A}^\vee}\ar[d] &
H_n^\vee\otimes\mathcal{O}_{\mathbb{P}^3}(1)\ar[r]\ar[d] & 0 & \\
& & 0 & 0 & &
}
\end{equation}

Thus $\mathbf{S}^0_n\subset MI_{n,n}$. In fact, the following result is true.

\begin{proposition}\label{Inn}
$\mathbf{S}^0_n=MI_{n,n}$. In particular, $MI_{n,n}$ is irreducible of dimension $3n^2+3n$, and hence $I_{n,n}$ is irreducible of dimension $2n^2+3n$.
\end{proposition}
\begin{proof}
We have to show that $MI_{n,n}\subset\mathbf{S}^0_n$. Let
$A\in MI_{n,n}$. Since $n=r$, by condition (i) from (\ref{space of nets})  the rank of the hyperweb of quadrics $A:H_n\otimes V\to H_n^\vee\otimes V^\vee$ is
$2n+2r=4n=\dim H_n^\vee\otimes V^\vee$, hence $A$ is invertible. By (\ref{S0n}), this means that
$A\in\mathbf{S}^0_n$.
\end{proof}

Now we proceed to spell out the relation between symplectic $(n,n)$-instantons and usual rank-2 instantons with second
Chern class $2n-1$. This relation is given at the level of spaces of hyperwebs of quadrics
$MI_{n,n}$ and $MI_{2n-1,1}$ interpreted as spaces of monads.

We need some more notation.
Let $B\in\mathbf{S}^0_n$. By definition, $B$ is an invertible anti-self-dual map
$H_n\otimes V\to H_n^\vee\otimes V^\vee$. Then the  inverse
\begin{equation}\label{B^-1}
B^{-1}:H_n^\vee\otimes V^\vee\to H_n\otimes V
\end{equation}
is also anti-self-dual.
Consider the vector space $\mathbf{\Sigma}_n=H_n^\vee\otimes H_{n-1}^\vee\otimes\wedge^2V^\vee$.
An element $C\in\mathbf{\Sigma}_n$ can be viewed as a linear map
$C:H_{n-1}\otimes V\to H_n^\vee\otimes V^\vee$,
and its transpose $C^\vee$ as a map
$C^\vee:H_n\otimes V\to H_{n-1}^\vee\otimes V^\vee$.
As the composition $C^\vee\circ B^{-1}\circ C$ is anti-self-dual, we can consider it as
an element of
$\wedge^2(H_{n-1}^\vee\otimes V^\vee)\simeq\mathbf{S}_{n-1}\oplus\wedge^2H_{n-1}^\vee\otimes S^2V^\vee$
(cf. (\ref{can decomp})).
Thus the condition
\begin{equation}\label{condition (i)}
C^\vee\circ B^{-1}\circ C\in\mathbf{S}_{n-1}
\end{equation}
makes sense.

Next, consider the upper horizontal triple in (\ref{diag E2n}) with $A=B$.
Twisting it by $\mathcal{O}_{\mathbb{P}^3}(1)$ and passing to global sections we obtain the exact triple
\begin{equation}\label{epsilon B}
0\to H_n\overset{{}^\sharp B}\to H_n^\vee\otimes\wedge^2V^\vee\overset{\epsilon(B)}\to H^0(E_{2n}(B)(1))\to 0
\end{equation}
Besides, interpreting $C\in\mathbf{\Sigma}_n$ as a map
${}^\sharp C:H_{n-1}\to H_n^\vee\otimes\wedge^2V^\vee$, we obtain the composition
$H_{n-1}\overset{{}^\sharp C}\to H_n^\vee\otimes\wedge^2V^\vee\overset{\epsilon(B)}\to H^0(E_{2n}(B)(1))$
which induces the morphism of sheaves
\begin{equation}\label{rhoA1A2}
\rho_{B,C}:\ H_{n-1}\otimes\mathcal{O}_{\mathbb{P}^3}(-1)\to E_{2n}(B).
\end{equation}
Note also that the maps $B:H_n\otimes V\to H_n^\vee\otimes V^\vee$
and
$C:H_{n-1}\otimes V\to H_n^\vee\otimes V^\vee$
provide a map
$(H_n\oplus H_{n-1})\otimes V\to H_n^\vee\otimes V^\vee$,
which induces the morphism of sheaves
\begin{equation}\label{tauBC}
\tau_{B,C}:\ (H_n\oplus H_{n-1})\otimes\mathcal{O}_{\mathbb{P}^3}(-1)\to
H_n^\vee\otimes V^\vee\otimes\mathcal{O}_{\mathbb{P}^3}.
\end{equation}

Now set
\begin{equation}\label{Xm}
X_n:=\left\{(B,C)\in\mathbf{S}_n^0\times\mathbf{\Sigma}_n\ \left|\
\begin{minipage}{20 em}
\begin{enumerate}
\item the\ condition\ (\ref{condition (i)})\ is\ satisfied,
\item $\rho_{B,C}$ in\ (\ref{rhoA1A2})\ is\ a\ subbundle\ inclusion,
\item $\tau_{B,C}$ in\ (\ref{tauBC})\ is\ a\ subbundle\ inclusion.
\end{enumerate}
\end{minipage}
\right.
\right\}
\end{equation}
By definition, $X_n$ is a locally closed subset of
$\mathbf{S}_n^0\times\mathbf{\Sigma}_n$.
Hence it is naturally endowed with a structure of a reduced scheme.

Now for any direct sum decomposition
\begin{equation}\label{xi}
\xi:H_{2n-1}\overset{\simeq}\to H_n\oplus H_{n-1},
\end{equation}
we obtain the corresponding decomposition
\begin{equation}\label{tilde xi}
\tilde{\xi}:\mathbf{S}_{2n-1}\overset{\simeq}\to\mathbf{S}_n\oplus\mathbf{\Sigma}_n\oplus\mathbf{S}_{n-1}:
A\mapsto (A_1(\xi),A_2(\xi),A_3(\xi)).
\end{equation}
Thus, considering the set $MI_{2n-1,1}$ of $(2n-1)$-instanton hyperwebs of quadrics as a subset of $\mathbf{S}_{2n-1}$,
we obtain a natural projection
\begin{equation}\label{prn fn0}
f_n:MI_{2n-1,1}\to\mathbf{S}_n\oplus\mathbf{\Sigma}_n:
A\mapsto (A_1(\xi),A_2(\xi)).
\end{equation}
The following result is proved in \cite[Theorems 1.1, 6.1 and Remark 7.6]{T}.
\begin{proposition}\label{isomorphism fn}
For a general decomposition $\xi$ in (\ref{xi}), there exists a dense open subset $MI_{2n-1,1}(\xi)$ of $MI_{2n-1,1}$
such that the projection $f_n$ in (\ref{prn fn0}) induces an isomorphism or integral schemes
\begin{equation}\label{isom fn}
f_n:MI_{2n-1,1}(\xi)\overset{\simeq}\to X_n:A\mapsto (A_1(\xi),A_2(\xi)).
\end{equation}
The inverse isomorphism is given by the formula
\begin{equation}\label{isom fn^-1}
f_n^{-1}:\ X_n\overset{\simeq}\to MI_{2n-1,1}(\xi):
\ (B,C)\mapsto\ \tilde{\xi}^{-1}(B,\ C,\ -C^\vee\circ B^{-1}\circ C).
\end{equation}
Besides, the projection
\begin{equation}\label{pr1 dom}
pr_1:\ X_n\to\mathbf{S}_n^0:\ (B,C)\mapsto B
\end{equation}
is dominant.
\end{proposition}

It is not hard to check that the morphism
$\rho_{B,C}:H_{n-1}\otimes\mathcal{O}_{\mathbb{P}^3}(-1)\to E_{2n}(B)$
defined in (\ref{rhoA1A2}) satisfies the condition ${}^t\rho_{B,C}\circ\rho_{B,C}=0$,
where ${}^t\rho_{B,C}$ is the composition
$$
{}^t\rho_{B,C}:E_{2n}(B)\overset{\phi}{\underset{\simeq}\to}
E_{2n}(B)^\vee\overset{\rho_{B,C}^\vee}\lra H_{n-1}^\vee\otimes\mathcal{O}_{\mathbb{P}^3}(1)
$$
and $\phi$ is a symplectic structure on $E_{2n}(B)$ (cf. \cite[formulas (71)-(72)]{T}). In other words,
we obtain an anti-self-dual monad
\begin{equation}\label{monad n}
0\to H_{n-1}\otimes\mathcal{O}_{\mathbb{P}^3}(-1)\overset{\rho_{B,C}}\lra E_{2n}(B)\overset{\phi}{\underset{\simeq}\to}
E_{2n}(B)^\vee\overset{\rho_{B,C}^\vee}\lra H_{n-1}^\vee\otimes\mathcal{O}_{\mathbb{P}^3}(1)\to0
\end{equation}
with cohomology sheaf
\begin{equation}\label{E2(B,C)}
E_2(A)=E_2(B,C):=\ker{}^t\rho_{B,C}/{\rm im\:}\rho_{B,C},\ \ \ A=f_n^{-1}(B,C).
\end{equation}
Next, by (\ref{pi nr}) we have the natural projection
\begin{equation}\label{pi}
\pi_{2n-1,1}: MI_{2n-1,1}\to I_{2n-1,1}:\ A\mapsto[E_2(A)].
\end{equation}
We have the following interpretation of the isomorphism (\ref{isom fn^-1}) on the
level of vector bundles:
\begin{equation}\label{E2 vers E2n}
[E_2(B,C)]=\pi_{2n-1,1}(f_n^{-1}(B,C)).
\end{equation}

\begin{remark}\label{quot monad}
Note that, according to the definitions (\ref{Monad A})-(\ref{space of nets}) of $MI_{2n-1,1}$ and $MI_{n,n}$, for any $A\in MI_{2n-1,1}$, if $B=A_1(\xi)$ is defined by the direct sum
decomposition (\ref{tilde xi}), one has two other anti-self-dual monads
\begin{equation}\label{again Monad A}
\mathcal{M}_A:\ \ 0\to H_{2n-1}\otimes\mathcal{O}_{\mathbb{P}^3}(-1)\overset{a_A}\to
W_{A}\otimes\mathcal{O}_{\mathbb{P}^3}\overset{a_A^\vee\circ q_A}
\to H_{2n-1}^\vee\otimes\mathcal{O}_{\mathbb{P}^3}(1)\to0
\end{equation}
\begin{equation}\label{Monad B}
\mathcal{M}_B:\ \ 0\to H_n\otimes\mathcal{O}_{\mathbb{P}^3}(-1)\overset{a_B}\to
W_B\otimes\mathcal{O}_{\mathbb{P}^3}\overset{a_B^\vee\circ q_B}
\to H_n^\vee\otimes\mathcal{O}_{\mathbb{P}^3}(1)\to0
\end{equation}
with cohomology sheaves
\begin{equation}\label{coho sheaves}
E_2(A)=\ker(a_A^\vee\circ q_A)/{\rm im\:} a_A,\
E_{2n}(B)=\ker(a_B^\vee\circ q_B)/{\rm im\:} a_B
\end{equation}
respectively. Moreover, \eqref{monad n} and \eqref{E2(B,C)} provide
an isomorphism
$w:\ W_B=H^2(E_2(B)\otimes\Omega_{\mathbb{P}^3})\overset{\simeq}\to H^2(E_{2n}(A)\otimes\Omega_{\mathbb{P}^3})=W_A$.
We thus obtain a commutative anti-self-dual diagram relating these monads:
\begin{equation}\label{two monads}
\xymatrix{
0\ar[r] & H_n\otimes\mathcal{O}_{\mathbb{P}^3}(-1)\ar[r]^-{a_B}\ar@{^{(}->}[d]^{i_\xi} &
W_B\otimes\mathcal{O}_{\mathbb{P}^3}\ar[r]^{q_B}_{\cong}\ar[d]^w_{\cong} &
W_B^\vee\otimes\mathcal{O}_{\mathbb{P}^3}\ar[r]^-{a_B^\vee} &
H_n^\vee\otimes\mathcal{O}_{\mathbb{P}^3}(1)\ar[r] & 0 \\
0\ar[r] & H_{2n-1}\otimes\mathcal{O}_{\mathbb{P}^3}(-1)\ar[r]^-{a_A} &
W_A\otimes\mathcal{O}_{\mathbb{P}^3}\ar[r]^{q_A}_{\cong} &
W_A^\vee\otimes\mathcal{O}_{\mathbb{P}^3}\ar[r]^-{a_A^\vee}\ar[u]^{w^\vee}_{\cong} &
H_{2n-1}^\vee\otimes\mathcal{O}_{\mathbb{P}^3}(1)\ar[r]\ar@{->>}[u]^{i_\xi^\vee} & 0,}
\end{equation}
where $i_\xi:H_n\hookrightarrow H_{2n-1}$ is the embedding induced by the decomposition (\ref{xi}).
In view of (\ref{coho sheaves}) and the canonical isomorphism $H_{2n-1}/i_\xi(H_n)\simeq H_{n-1}$,  from this diagram
we obtain the monad
\begin{equation}\label{quotient monad}
\mathcal{M}_{A,B}:\ \ 0\to H_{n-1}\otimes\mathcal{O}_{\mathbb{P}^3}(-1)\overset{a_{A,B}}\to
E_{2n}(B)\overset{\phi}{\underset{\simeq}\to}E_{2n}(B)^\vee
\overset{a_{A,B}^\vee}\to H_{2n-1}^\vee\otimes\mathcal{O}_{\mathbb{P}^3}(1)\to0
\end{equation}
with cohomology sheaf
\begin{equation}\label{rk2 coho sheaf}
E_2(A)=\ker(a_{A,B}^\vee\circ \phi)/{\rm im\:} a_A.
\end{equation}
We call (\ref{quotient monad}) the {\it quotient monad} of the monads (\ref{again Monad A}) and (\ref{Monad B}).
\end{remark}
\begin{remark}\label{parameters for diag}
Note that, by Proposition \ref{isomorphism fn}, the set of all diagrams (\ref{two monads}) is parametrized
by the irreducible variety $I_{2n-1,1}(\xi)$.
\end{remark}

\end{sub}

\begin{sub}{\bf Example: a special family of symplectic $(n,r)$-instantons.}
\label{nr-instantons}
\rm
Now assume $n\ge2$ and, for any integer $r,\ 2\le r\le n-1$, consider an inclusion
\begin{equation}\label{mono tau}
\tau:H_{2n-r}\hookrightarrow H_{2n-1}
\end{equation}
such that
\begin{equation}\label{tau contains ixi}
\tau(H_{2n-r})\supset i_\xi(H_n).
\end{equation}
We obtain a hyperweb of quadrics
$$
A_\tau\in S^2H_{2n-r}^\vee\otimes\wedge^2V^\vee
$$
as the image of $A$ under the map
$S^2H_{2n-1}^\vee\otimes\wedge^2V^\vee \to S^2H_{2n-r}^\vee\otimes\wedge^2V^\vee$
induced by $\tau$.
The corresponding monad
\begin{equation}\label{3rd monad}
\mathcal{M}_\tau:\ \ 0\to H_{2n-r}\otimes\mathcal{O}_{\mathbb{P}^3}(-1)\overset{a_\tau}\to
W_{A}\otimes\mathcal{O}_{\mathbb{P}^3}\overset{a_\tau^\vee\circ q_A}
\to H_{2n-r}^\vee\otimes\mathcal{O}_{\mathbb{P}^3}(1)\to0,
\end{equation}
has a rank-$2r$ cohomology bundle
\begin{equation}\label{rk2n coho sheaf}
E_{2r}(A_\tau)=\ker(a_\tau^\vee\circ q_A)/{\rm im\:} a_\tau.
\end{equation}
where $a_\tau:=a_A\circ\tau$. By construction, $E_{2r}(A_\tau)$ inherits a natural symplectic structure
\begin{equation}\label{sympl E2r}
\phi_r:\ E_{2r}(A_\tau)\overset{\simeq}\to E_{2r}(A_\tau)^\vee.
\end{equation}
Besides, in view of (\ref{tau contains ixi}), the monad (\ref{3rd monad}) can be inserted as a midle row into the diagram (\ref{two monads}), extending it to a three-row
commutative anti-self-dual diagram. Arguing
as in Remark \ref{quot monad} we obtain, in addition to the quotient monad (\ref{quotient monad}), two more quotient
monads:
\begin{equation}\label{2nd quotient monad}
\mathcal{M}'_{\tau}:\ \ 0\to H_{n-r}\otimes\mathcal{O}_{\mathbb{P}^3}(-1)\overset{a'_\tau}\to
E_{2n}(B)\overset{\phi}{\underset{\simeq}\to}E_{2n}(B)^\vee
\overset{{a'}_{\tau}^\vee}\to H_{n-r}^\vee\otimes\mathcal{O}_{\mathbb{P}^3}(1)\to0,
\end{equation}
$$
E_{2r}(A_\tau)=\ker({a'}_{\tau}^\vee\circ \phi)/{\rm im\:} a'_\tau,
$$
\begin{equation}\label{3rd quotient monad}
\mathcal{M}''_{\tau}:\ \ 0\to H_{r-1}\otimes\mathcal{O}_{\mathbb{P}^3}(-1)\overset{a''_\tau}\to
E_{2r}(B)\overset{\phi_\tau}{\underset{\simeq}\to}E_{2r}(B)^\vee
\overset{{a''}_{\tau}^\vee}\to H_{r-1}^\vee\otimes\mathcal{O}_{\mathbb{P}^3}(1)\to0,
\end{equation}
$$
E_2(A)=\ker({a'}_{\tau}^\vee\circ \phi_\tau)/{\rm im\:} a_A.
$$
From (\ref{vanish cohom of E2n}) and (\ref{2nd quotient monad}) we easily deduce:
\begin{equation}\label{vanish cohom of E2r}
h^0(E_{2r}(A_\tau))=h^i(E_{2r}(A_\tau)(-2))=0,\ \ \ i\ge0,\ \ \ c_2(E_{2r}(A_\tau))=2n-r.
\end{equation}
By definition, this together with (\ref{3rd monad})-(\ref{sympl E2r}) means that
\begin{equation}\label{belongs to E2r}
[E_{2r}(A_\tau)]\in I_{2n-r,r}.
\end{equation}

\begin{remark}\label{Dnr}
Observe that, in view of (\ref{mono tau}), the maps $\tau$ belong to the set
$$
N_{n,r}:=\{\tau\in{\rm Hom}(H_{2n-r},H_{2n-1})|\ \tau\ {\rm\ is\ injective\ and}\ {\rm im\:}\tau
\supset{\rm im\:}i_\xi\}.
$$
When $A\in MI_{2n-1,1}(\xi)$ is fixed,$N_{n,r}$ parametrizes some family of hyperwebs $A_\tau$ from $MI_{2n-r,r}$.
Since $N_{n,r}$ is a principal $GL(H_{2n-r})$-bundle over an open subset of the Grassmannian $Gr(n-r,n-1)$, it it is irreducible. Thus, by Remark \ref{parameters for diag}, the family of the three-row extensions of the diagrams (\ref{two monads}) can be parametrized by the irreducible variety
$MI_{2n-1,1}(\xi)\times N_{n,r}$.
Hence the family $D_{n,r}$ of isomorphism classes of symplectic rank-$2r$ bundles obtained from these
diagrams by formula (\ref{rk2n coho sheaf}) is an irreducible locally closed subset of $I_{2n-r,r}$.

Note that it is a priori not clear whether the closure of $D_{n,r}$ in $I_{2n-r,r}$ is an irreducible component of
$I_{2n-r,r}$.
\end{remark}

\begin{definition}\label{property (*)}
Let $2\le r\le n-1$. We say that {\it $A\in MI_{2n-r,r}$ satisfies property} (*) if there exists a
monomorphism $i:H_n\hookrightarrow H_{2n-r}$ such that the image $B$ of $A$ under the surjection
$\mathbf{S}_{2n-r}\twoheadrightarrow\mathbf{S}_n$ induced by $i$ is invertible as a homomorphism
$B:H_n\otimes V\to H_n^\vee\otimes V^\vee$.

The property (*) is clearly an open condition on $A$.
Moreover, since $\pi_{2n-r,r}: MI_{2n-r,r}\to I_{2n-r,r}$ is a principal bundle (Theorem \ref{princbundle}),
if an element $A\in\pi_{2n-r,r}^{-1}([E_{2r}])$ satisfies (*), then any
other point $A'\in\pi_{2n-r,r}^{-1}([E_{2r}])$ satisfies (*).
We thus say that {\it $[E_{2r}]\in I_{2n-r,r}$ satisfies property} (*) if some (hence any)
$A\in\pi_{2n-r,r}^{-1}([E_{2r}])$ satisfies property (*). It is obviously an open condition on
$[E_{2r}]\in I_{2n-r,r}$.
\end{definition}

\begin{remark}\label{I*}
By Proposition \ref{isomorphism fn} and using (\ref{tau contains ixi}), we see that any
$[E_{2r}]\in D_{n,r}$, as well as any $A\in f_n^{-1}(D_{n,r})$ satisfies property (*). We define
\begin{equation}\label{def of I*}
I_{2n-r,r}^*:=I_{(1)}\cup\ldots\cup I_{(k)},
\end{equation}
where $I_{(1)},\ldots, I_{(k)}$
are all the irreducible components of $I_{2n-r,r}$ whose general points
satisfy property (*). By definition, $D_{n,r}\subset I_{2n-r,r}^*$, hence $I_{2n-r,r}^*$
is nonempty. We also set $MI_{2n-r,r}^*= \pi_{2n-r,r}^{-1}(I_{2n-r,r}^*)$, so that the map
$\pi_{2n-r,r}:MI_{2n-r,r}^*\to I_{2n-r,r}^*$ is a principal bundle with structure group
$GL(H_{2n-r})/\{\pm 1\}$.

\end{remark}

\end{sub}

\bigskip

\section{Irreducibility of $I^*_{2n-r,r}$}\label{I^*nr}

\begin{sub}{\bf A dense open subset $X_{n,r}$ of $MI_{2n-r,r}^*$. Reduction of the irreducibility
of $I^*_{n,r}$ to that of $X_{n,r}$.}\label{I^*nr and Xnr}\rm\
In this section we prove the irreducibility of the component $I_{2n-r,r}^*$ of $I_{2n-r,r}$ defined in (\ref{def of I*}),
see Theorem \ref{Irred of I*}.
The explicit construction of symplectic instantons in Section \ref{explicit} gives us a hint to the proof. We proceed along the lines of Subsection \ref{nn-instantons}.

Take any $B\in\mathbf{S}^0_n$ and consider it as an invertible anti-self-dual
linear map $H_n\otimes V\to H_n^\vee\otimes V^\vee$. Then $B^{-1}$ is also anti-self-dual.
Let
\begin{equation}\label{Sigma nr}
\mathbf{\Sigma}_{n,r}:=H_{n-r}^\vee\otimes H_n^\vee\otimes\wedge^2 V^\vee.
\end{equation}
An element $C\in\mathbf{\Sigma}_n$ can be understood as a map
$C:H_{n-r}\otimes V\to H_n^\vee\otimes V^\vee$,
and its transpose $C^\vee$ is a map
$H_n\otimes V\to H_{n-r}^\vee\otimes V^\vee$.
The composition $C^\vee\circ B^{-1}\circ C$ is anti-self-dual, i.e.,
it is an element of
$\wedge^2(H_{n-r}^\vee\otimes V^\vee)\simeq\mathbf{S}_{n-r}\oplus\wedge^2H_{n-r}^\vee\otimes S^2V^\vee$
(cf. (\ref{can decomp})).
We will later impose the condition
\begin{equation}\label{condition (ir)}
C^\vee\circ B^{-1}\circ C\in\mathbf{S}_{n-r}.
\end{equation}

Next, as in (\ref{epsilon B}), we have a well defined epimorphism
$\epsilon(B):H_n^\vee\otimes\wedge^2V^\vee\twoheadrightarrow H^0(E_{2n}(B)(1))$.
Besides, interpreting the above element  $C\in\mathbf{\Sigma}_{n,r}$ as a map
${}^\sharp C:H_{n-r}\to H_n^\vee\otimes\wedge^2V^\vee$, we obtain the composition
$H_{n-r}\overset{{}^\sharp C}\to H_n^\vee\otimes\wedge^2V^\vee\overset{\epsilon(B)}\to H^0(E_{2n}(B)(1))$
which induces the morphism of sheaves
\begin{equation}\label{rhoBC}
\rho_{B,C}:\ H_{n-r}\otimes\mathcal{O}_{\mathbb{P}^3}(-1)\to E_{2n}(B).
\end{equation}
Note also that $B:H_n\otimes V\to H_n^\vee\otimes V^\vee$
and
$C:H_{n-r}\otimes V\to H_n^\vee\otimes V^\vee$
define a map
$(H_n\oplus H_{n-r})\otimes V\to H_n^\vee\otimes V^\vee$
which induces the morphism of sheaves
\begin{equation}\label{tauA1A2}
\tau_{B,C}:\ (H_n\oplus H_{n-r})\otimes\mathcal{O}_{\mathbb{P}^3}(-1)\to
H_n^\vee\otimes V^\vee\otimes\mathcal{O}_{\mathbb{P}^3}.
\end{equation}

Now set
\begin{equation}\label{Xm1}
X_{n,r}:=\left\{(B,C)\in\mathbf{S}_n^0\times\mathbf{\Sigma}_{n,r}\ \left|\
\begin{minipage}{20 em}
\begin{enumerate}
\item the\ condition\ (\ref{condition (ir)})\ is\ satisfied,
\item $\rho_{B,C}$ in\ (\ref{rhoBC})\ is\ a\ subbundle\ inclusion,
\item $\tau_{B,C}$ in\ (\ref{tauA1A2})\ is\ a\ subbundle\ inclusion.
\end{enumerate}
\end{minipage}
\right.
\right\}
\end{equation}
By definition, $X_{n,r}$ is a locally closed subset of
$\mathbf{S}_n^0\times\mathbf{\Sigma}_{n,r}$.
Hence it has a natural structure of  reduced scheme.

Now for an arbitrary direct sum decomposition
\begin{equation}\label{new xi}
\xi:H_{2n-r}\overset{\simeq}\to H_n\oplus H_{n-r}
\end{equation}
we obtain the corresponding decomposition
\begin{equation}\label{new tilde xi}
\tilde{\xi}:\mathbf{S}_{2n-r}\overset{\simeq}\to\mathbf{S}_n\oplus\mathbf{\Sigma}_{n,r}\oplus\mathbf{S}_{n-r}:
A\mapsto (A_1(\xi),A_2(\xi),A_3(\xi)).
\end{equation}
Thus, considering the set $MI_{2n-r,r}$ of symplectic $(2n-r,r)$-instanton hyperwebs of quadrics as a subset of
$\mathbf{S}_{2n-r}$,
we obtain a natural projection
\begin{equation}\label{prn fn}
f_{n,r}:MI_{2n-r,r}\to\mathbf{S}_n\oplus\mathbf{\Sigma}_{n,r}:
A\mapsto (A_1(\xi),A_2(\xi)).
\end{equation}
We now prove the following result parallel to Proposition \ref{isomorphism fn}.
\begin{theorem}\label{Xnr isom MI2n-r}
Let $n\ge3$ and $2\le r\le n-1$.

(i) For a general decomposition $\xi$ in (\ref{new xi}) there is an open dense subset
$MI_{2n-r,r}^*(\xi)$ of $MI_{2n-r,r}^*$ and an isomorphism of reduced schemes
\begin{equation}\label{isomorphism fnr}
f_{n,r}:\ MI_{2n-r,r}^*(\xi)\overset{\simeq}\to X_{n,r}:
\ A\mapsto(A_1(\xi),A_2(\xi)),
\end{equation}
where $A_1(\xi)$ and $A_2(\xi)$ are defined by (\ref{new tilde xi}).

(ii) The inverse isomorphism is given by the formula
\begin{equation}\label{isomorphism fnr^-1}
f_{n,r}^{-1}:\ X_{n,r}\overset{\simeq}\to MI_{2n-r,r}^*(\xi):
\ (B,C)\mapsto\ \widetilde{\xi}^{-1}(B,\ C,\ -C^\vee\circ B^{-1}\circ C),
\end{equation}
where $\widetilde{\xi}$ is defined by (\ref{new tilde xi}).
\end{theorem}

\begin{proof}
Set $MI_{2n-r,r}^*(\xi):=\{A\in MI_{2n-r,r}^*\ |\ A$ satisfies property (*) for the monomorphism
$i:H_n\hookrightarrow H_{2n-r}$ defined by $\xi\}$. It follows from
Definition \ref{property (*)}and Remark \ref{I*} that, for a general decomposition $\xi$ in (\ref{new xi}),
$MI_{2n-r,r}^*(\xi)$ is a dense open subset of $MI_{2n-r,r}^*$. Then, for this choice of $\xi$, the proof of this
Theorem essentially mimics the proof of \cite[Proposition 6.1]{T} in which we make the  substitution
$m+1\mapsto n,\ m\mapsto n-r$ and change the notation accordingly.
\end{proof}

The proof of the following theorem will be given in Subsection \ref{Irreducibility of Xnr}.

\begin{theorem}\label{Irred of Xnr}
$X_{n,r}$ is irreducible of dimension $(2n-r)^2+4(2n-r)(r+1)-r(2r+1)$.
\end{theorem}

From Theorems \ref{Xnr isom MI2n-r} and \ref{Irred of Xnr} it follows that $MI_{2n-r,r}^*$ is irreducible
of dimension $(2n-r)^2+4(2n-r)(r+1)-r(2r+1)$ for any $n\le3$ and $2\le r\le n-1$.
Hence $I_{2n-r,r}^*$ is irreducible
of dimension $4(2n-r)(r+1)-r(2r+1)$ for these values of $n$ and $r$. Note that the irreducibility of $I_{2n-r,r}^*$ is
also true when $r=n$, and in this case $I_{n,n}^*$ coincides with $I_{n,n}$.  Substituting
$2n-1\mapsto n$, we obtain the following main result of the paper.

\begin{theorem}\label{Irred of I*}
For any integer $r\ge2$ and for any integer $n\ge r$ such that $n\equiv r({\rm mod}2)$, $I_{n,r}^*$
is an irreducible component of $I_{n,r}$ of dimension $4n(r+1)-r(2r+1)$.
\end{theorem}

\end{sub}

\medskip

\begin{sub}{\bf Proof of the irreducibility of $X_{n,r}$.}\label{Irreducibility of Xnr}
\rm\ In this subsection we give the proof of Theorem \ref{Irred of Xnr}. Define
\begin{equation}\label{tilde Xm}
\widetilde{X}_{n,r}:=
\{(D,C)\in (\mathbf{S}_n^\vee)^0\times\mathbf{\Sigma}_{n,r}\ |\
\ (C^\vee\circ D\circ C:H_{n-r}\otimes V\to H_{n-r}^\vee\otimes V^\vee)
\in\mathbf{S}_{n-r}\},
\end{equation}
a closed subscheme of
$(\mathbf{S}_m^\vee)^0\times\mathbf{\Sigma}_{n,r}$
defined by the equations
\begin{equation}\label{eqns of tilde Xm}
C^\vee\circ D\circ C\in \mathbf{S}_{n-r}.
\end{equation}
Since the conditions (ii) and (iii) in the definition (\ref{Xm}) of $X_{n,r}$ are open and $X_{n,r}$
is nonempty (see Theorem \ref{Xnr isom MI2n-r}), the isomorphism
$$
\mathbf{S}_n^0\overset{\simeq}\to(\mathbf{S}_n^\vee)^0:\ B\mapsto B^{-1}
$$
implies that $X_{n,r}$ is a
nonempty open subset of $(\widetilde{X}_{n,r})_{red}$,
\begin{equation}\label{open in tilde Xm}
\emptyset\ne X_{n,r}\intoo{\rm open}(\widetilde{X}_{n,r})_{red}.
\end{equation}

Fix a direct sum decomposition
$$
H_n\overset{\simeq}\to H_{n-r}\oplus H_r.
$$
Then any linear map
\begin{equation}\label{map C}
C\in\mathbf{\Sigma}_{n,r}=\Hom(H_{n-r},H_n^\vee\otimes\wedge^2 V^\vee),\ \ \
C:H_{n-r}\otimes V\to H_n^\vee\otimes V^\vee,\ \ \ \
\end{equation}
can be represented as a map
\begin{equation}\label{decompn C}
C:H_{n-r}\otimes V\to H_{n-r}^\vee\otimes V^\vee\ \oplus\ H_r^\vee\otimes V^\vee,
\end{equation}
or else as a block matrix
\begin{equation}\label{matrix of C}
C= \left(
\begin{array}{c}
\phi \\
\psi
\end{array}
\right),
\end{equation}
where
\begin{equation}\label{phi, b}
\phi\in\Hom(H_{n-r},H_{n-r}^\vee)\otimes\wedge^2 V^\vee=\mathbf{\Phi}_{n-r},\ \ \
\psi\in\mathbf{\Psi}_{n,r}:=
\Hom(H_{n-r},H_r^\vee)\otimes\wedge^2 V^\vee.
\end{equation}
Similarly, any
$D\in (\mathbf{S}^\vee_n)^0\subset\mathbf{S}^\vee_n= S^2H_n\otimes\wedge^2V\subset
\Hom(H_n^\vee\otimes V^\vee,H_n\otimes V)$
can be represented in the form

\begin{equation}\label{matrix of D}
D= \left(
\begin{array}{cc}
D_1 & \lambda \\
-\lambda^\vee & \mu
\end{array}
\right),
\end{equation}
where
\begin{equation}\label{A, lambda, mu}
D_1\in\mathbf{S}^\vee_{n-r}\subset
\Hom(H_{n-r}^\vee\otimes V^\vee,H_{n-r}\otimes V),\ \ \
\end{equation}
$$
\lambda\in\mathbf{L}_{n,r}:=\Hom(H_r^\vee,H_{n-r})\otimes\wedge^2 V, \ \ \
\mu\in\mathbf{M}_r:=S^2H_r\otimes\wedge^2 V.
$$
By (\ref{matrix of C}) and (\ref{matrix of D}) the composition
$$
C^\vee\circ D\circ C:H_{n-r}\otimes V\to
H_{n-r}^\vee\otimes V^\vee\ \
(C^\vee\circ D\circ C\in\wedge^2(H_{n-r}^\vee\otimes V^\vee))
$$
can be written in the form
\begin{equation}\label{CDC}
C^\vee\circ D\circ C=\phi^\vee\circ D_1\circ\phi+\phi^\vee\circ\lambda\circ\psi-
\psi^\vee\circ\lambda^\vee\circ\phi+\psi^\vee\circ\mu\circ\psi.
\end{equation}
By (\ref{matrix of C})-(\ref{A, lambda, mu}) we have
$$
\mathbf{S}^\vee_n\times\mathbf{\Sigma}_{n,r}=\mathbf{S}^\vee_{n-r}\times\mathbf{\Phi}_{n-r}
\times\mathbf{\Psi}_{n,r}\times\mathbf{L}_{n,r}\times\mathbf{M}_r,
$$
and there are well defined morphisms
$$
\tilde{p}:\widetilde{X}_{n,r}\to\mathbf{L}_{n,r}\times\mathbf{M}_r:(D_1,\phi,\psi,\lambda,\mu)\mapsto(\lambda,\mu).
$$
and
$$
p:=\tilde{p}|\overline{X}_{n,r}:\overline{X}_{n,r}\to\mathbf{L}_{n,r}\oplus\mathbf{M}_r,
$$
where $\overline{X}_{n,r}$ is the closure of $X_{n,r}$ in $(\mathbf{S}_n^\vee)^0\times\mathbf{\Sigma}_{n,r}$.
We now invoke the following result from \cite{T}:

\begin{proposition}\label{nondeg for general}
Let $n\ge2$. Then for any $D\in(\mathbf{S}^\vee_n)^0$ and for a general choice of the decomposition
$H_n\overset{\sim}\to H_{n-r}\oplus H_r$, the block
$D_1$ of $D$ in (\ref{matrix of D}) is nondegenerate.
\end{proposition}

\begin{proof} See \cite[Proposition 7.3]{T}.
By repeatedly applying this proposition $r$ times, we can find a decomposition
$H_n\overset{\sim}\to H_{n-r}\oplus H_r$ such that
$D_1:H_{n-r}^\vee\otimes V^\vee\to H_{n-r}\otimes V$ in (\ref{matrix of D})
is nondegenerate, i.e., $D_1\in(\mathbf{S}^\vee_{n-r})^0$.
\end{proof}

Let $\mathcal{X}$ be any irreducible component of $X_{n,r}$ and let $\overline{\mathcal{X}}$ be its closure in
$\overline{X}_{n,r}$. Fix a point
$z=(D_1,\phi,\psi,\lambda,\mu)\in \mathcal{X}$ not lying in the components of $X_{n,r}$ different from $\mathcal{X}$.
Consider the morphism
\begin{equation}\label{f}
f:\ \mathbb{A}^1\to\overline{\mathcal{X}}:\ t\mapsto(D_1,t^2\phi,t\psi,t\lambda,t^2\mu),\ \ \ f(1)=z,
\end{equation}
which is well defined by (\ref{CDC}). By definition, the point
$f(0)=(D_1,0,0,0,0)$ lies in the fibre $p^{-1}(0,0)$. Hence,
$p^{-1}(0,0)\cap\overline{\mathcal{X}}\ne\emptyset$.
In other words,
\begin{equation}\label{nonempty fibre}
\rho^{-1}(0,0)\ne\emptyset,\ \ \ \ \text{\rm where}\ \ \ \rho:=p|\overline{\mathcal{X}}.
\end{equation}
Now, it follows from (\ref{CDC}) and the definition of $\widetilde{X}_{n,r}$ that
\begin{equation}\label{zero fibre}
\tilde{p}^{-1}(0,0)=\{(D_1,\phi,\psi)\in(\mathbf{S}^\vee_{n-r})^0\times\mathbf{\Phi}_{n-r}\times\mathbf{\Psi}_{n,r}\ |\
\phi^\vee\circ D_1\circ\phi\in\mathbf{S}_{n-r}\}.
\end{equation}
Consider the set
$$
Z_{n-r}=\{(D,\phi)\in(\mathbf{S}^\vee_{n-r})^0\times\mathbf{\Phi}_{n-r}\ |\
\phi^\vee\circ D\circ\phi\in\mathbf{S}_{n-r}\}.
$$
It carries a natural scheme structure, where  it is a closed subscheme of
$(\mathbf{S}^\vee_{n-r})^0\times\mathbf{\Phi}_{n-r}$.
Comparing the definition of $Z_{n-r}$ with (\ref{zero fibre}) we see that there are scheme-theoretic inclusions of schemes
\begin{equation}\label{zero fibre subset Zm times Psim}
\rho^{-1}(0,0)\subset p^{-1}(0,0)\subset
\tilde{p}^{-1}(0,0)=Z_{n-r}\times\mathbf{\Psi}_{n,r}.
\end{equation}
By \cite[Theorem 7.2]{T}, $Z_{n-r}$ is an integral scheme of dimension $4(n-r)(n-r+2)$. This together with
(\ref{zero fibre subset Zm times Psim}) implies that
\begin{equation}\label{dim fibre le dim Zm+dim Psim}
\dim\rho^{-1}(0,0)\le\dim p^{-1}(0,0)\le\dim Z_{n-r}+\dim\mathbf{\Psi}_{n,r}=4(n-r)(n-r+2)+6r(n-r)=
\end{equation}
$$
=(n-r)(4n+2r+8).
$$
Hence in view of (\ref{nonempty fibre})
\begin{equation}\label{dim X le}
\dim \overline{\mathcal{X}}\le\dim\rho^{-1}(0,0)+\dim\mathbf{L}_{n,r}+\dim\mathbf{M}_r\le(n-r)(4n+2r+8)+6r(n-r)+3r(r+1)=
\end{equation}
$$
=(2n-r)^2+4(2n-r)(r+1)-r(2r+1).
$$
On the other hand, formula (\ref{dim MInr ge...}), with $2n-r$ substituted for $n$, and
Theorem \ref{Xnr isom MI2n-r}(ii) show that, for any point $x\in\mathcal{X}$ such that
$A:=f_{n,r}^{-1}(x)\in MI_{2n-r,r}^0(\xi)$,
\begin{equation}\label{dim X ge}
(2n-r)^2+4(2n-r)(r+1)-r(2r+1)\le\dim_A MI_{2n-r,r}^0(\xi)=\dim\overline{\mathcal{X}}.
\end{equation}
Comparing (\ref{dim X le}) with (\ref{dim X ge}), we see that all the inequalities in
(\ref{dim fibre le dim Zm+dim Psim})-(\ref{dim X ge}) are equalities. In particular,
\begin{equation}\label{dim fibre=dim X'-dim base}
\dim\rho^{-1}(0,0)=\dim(Z_{n-r}\times\mathbf{\Psi}_{n,r})=
\dim\overline{\mathcal{X}}-\dim(\mathbf{L}_{n,r}\times\mathbf{M}_r).
\end{equation}
Since by Theorem \cite[Theorem 7.2]{T} the scheme $Z_{n-r}$ is integral and so $Z_{n-r}\times\mathbf{\Psi}_{n,r}$ is
integral as well, (\ref{zero fibre subset Zm times Psim}) and (\ref{dim fibre=dim X'-dim base}) yield the equalities
of integral schemes
\begin{equation}\label{zero fibre = Zm times Psim}
\rho^{-1}(0,0)=p^{-1}(0,0)=
\tilde{p}^{-1}(0,0)=Z_{n-r}\times\mathbf{\Psi}_{n,r}.
\end{equation}

Now we invoke one auxiliary result from \cite{T}.
\begin{lemma}\label{flat implies irred}
Let $f:X\to Y$ be a morphism of reduced schemes, where $Y$ is a smooth integral scheme. Assume that there exists a
closed point $y\in Y$ such that for any irreducible component $X'$ of $X$ the following conditions are satisfied:

(a) $\dim f^{-1}(y)=\dim X'-\dim Y$,

(b) the scheme-theoretic inclusion of fibres $(f|_{X'})^{-1}(y)\subset f^{-1}(y)$ is an isomorphism of integral
schemes.\\
Then

(i) there exists an open subset $U$ of\ $Y$ containing the point $y$ such that the morphism
$f|_{f^{-1}(U)}:f^{-1}(U)\to U$ is flat, and

(ii) $X$ is integral.
\end{lemma}

\begin{proof}
See \cite[Lemma 7.4]{T}.
\end{proof}

Applying assertions (i)-(ii) of this lemma to
$X=X_{n,r},\ X'=\mathcal{X},\ Y=\mathbf{L}_{n,r}\times\mathbf{M}_r,\ y=(0,0),f=p$, and using
(\ref{dim fibre=dim X'-dim base}) and (\ref{zero fibre = Zm times Psim}), we obtain that
$X_{n,r}$ is integral of dimension $(2n-r)^2+4(2n-r)(r+1)-r(2r+1)$. Theorem \ref{Irred of Xnr} is proved.

\end{sub}

\bigskip

{

}

\end{document}